\documentclass[11pt,reqno,twoside]{amsart}
\usepackage[asymmetric,top=2.5cm,bottom=2.8cm,left=2.38cm,right=2.38cm]{geometry}
\usepackage{amsmath,amssymb,amsfonts,amsthm,mathrsfs,bm}

\newtheorem{theorem}{Theorem}[section]
\newtheorem{lemma}[theorem]{Lemma}
\newtheorem{proposition}[theorem]{Proposition}

\theoremstyle{definition}
\newtheorem{definition}[theorem]{Definition}
\newtheorem{remark}[theorem]{Remark}

\DeclareMathOperator{\curl}{curl}
\DeclareMathOperator{\diver}{div}

\newcommand{\0}{\mathbf{0}}
\newcommand{\R}{\mathbb{R}}
\newcommand{\C}{\mathbb{C}}
\newcommand{\Sph}{\mathbb{S}^2}
\newcommand{\E}{\mathbf{E}}
\newcommand{\Hf}{\mathbf{H}}
\newcommand{\x}{\mathbf{x}}

\usepackage{graphicx}
\usepackage{amssymb}
\usepackage{epstopdf}
\usepackage{float}
\usepackage{caption}
\usepackage{mathtools}
\usepackage{booktabs}
\usepackage{subcaption}
\usepackage{amsmath, amssymb, amsthm}
\usepackage{mathtools}
\usepackage{geometry}
\usepackage{hyperref}
\hypersetup{colorlinks=true, linkcolor=blue, citecolor=blue, urlcolor=blue}

\usepackage{graphicx}
\usepackage{subcaption}

\usepackage[T1]{fontenc}
\usepackage{lmodern}

\usepackage{lipsum} 
\usepackage{hyperref}
\hypersetup{hypertex=true,
	colorlinks=true,
	linkcolor=blue,
	anchorcolor=blue,
	citecolor=blue}

\usepackage[usenames]{color}

\usepackage{xcolor}

\theoremstyle{definition}

\numberwithin{equation}{section}

\allowdisplaybreaks

\usepackage{xstring} 

\newcommand{\AuthorInfo}[4]{%
	\textsc{#1}%
	\IfStrEq{#4}{true}{$^{*}$}{}\\
	#2\\
	\textit{E-mail:} \texttt{#3}%
	\IfStrEq{#4}{true}{\\\textit{$^{*}$Corresponding author.}}{}
	\par\vspace{0.5em}
}

\title[Localized gradient enhancement near anisotropic electromagnetic scatterers]
{Localized gradient enhancement near anisotropic electromagnetic scatterers}

\date{} 

\begin{document}

	\maketitle
	
	\begin{center}
		\bigskip
		\footnotesize
		\AuthorInfo{Weisheng Zhou}{School of Mathematics, Jilin University, Changchun 130012, People's Republic of China}{wszhou1211@163.com, zhouws24@mails.jlu.edu.cn}{false}

		\AuthorInfo{Huaian Diao}{School of Mathematics and Key Laboratory of Symbolic Computation and Knowledge Engineering of Ministry of Education, Jilin University, Changchun 130012, People's Republic of China}{diao@jlu.edu.cn, hadiao@gmail.com}{true}
		
		\AuthorInfo{Hongyu Liu}{Department of Mathematics, City University of Hong Kong, Kowloon, Hong Kong SAR, People's Republic of China}{hongyu.liuip@gmail.com, hongyliu@cityu.edu.hk}{false}
	\end{center}
	\normalsize
	
	\begin{abstract}
		This work investigates time-harmonic electromagnetic scattering governed by the Maxwell system, where bounded anisotropic scatterers are embedded in a homogeneous electromagnetic background. We focus on the localized enhancement of the gradients of the total electric and magnetic fields in small boundary-attached neighborhoods of finitely many prescribed points near boundaries of anisotropic electromagnetic scatterers. We show that, through a suitable construction of incident electromagnetic waves, the gradients of both the total electric field and the total magnetic field can be made arbitrarily large in these neighborhoods.  The main strategy is based on the introduction of auxiliary boundary-attached electromagnetic neighborhoods and the associated electric and magnetic fields, which exhibit strong gradient variation near the prescribed points. Using the approximation property of Maxwell Herglotz wave functions, these auxiliary fields are then approximated by physically admissible incident waves in the neighborhood of the scatterers. Together with the well-posedness and continuous dependence of the anisotropic scattering problem, this implies that the corresponding scattered field can be controlled to be sufficiently weak in the relevant  region. Consequently, the total field is dominated by the incident field near the prescribed points and inherits its large-gradient behavior. The result provides a theoretical mechanism for localized gradient enhancement in anisotropic electromagnetic scattering and may have implications for field concentration, high-resolution probing, and sensitivity analysis of electromagnetic responses in complex media.
		
		
	\end{abstract}
	\maketitle

	\section{Introduction}
	
	Localized field concentration is a fundamental mechanism in wave propagation,
	electromagnetic engineering, inverse scattering, and high-resolution wave
	imaging. In engineering applications, localized electromagnetic enhancement is
	commonly pursued through near-field focusing metasurfaces, programmable
	metasurfaces, and dynamic metasurface antennas, where the phase, amplitude, or
	polarization response of structured media is designed to concentrate
	electromagnetic energy in prescribed near-field regions
	\cite{WuZhangCuiZhangXu2024,DingEtAl2025,ZhangShlezingerGuidiDardariImaniEldar2025}.
	Related ideas also appear in optical metasurfaces and structured-light
	technologies, where engineered wavefronts are used to tailor light--matter
	interaction and improve high-resolution imaging or local response
	\cite{ChenWenQiu2019,MaEtAl2024}. In electromagnetic media with complex material
	parameters, it is therefore important not only to enhance the field amplitude,
	but also to generate strong local spatial variation of the field near prescribed
	locations. Such a gradient-level enhancement is closely related to
	high-sensitivity probing, enhanced material response, and the controllability of
	wave--matter interaction in structured electromagnetic environments.

	In this paper, we study a time-harmonic Maxwell scattering problem for bounded
	anisotropic electromagnetic scatterers embedded in  a homogeneous exterior background. The scatterer may
	consist of finitely many pairwise disjoint components, and we prescribe finitely
	many target points on the outer boundary of the scatterer.
	Our objective is to construct physically admissible incident electromagnetic
	waves such that the gradients of the corresponding total electric and magnetic
	fields become arbitrarily large in small boundary-attached neighborhoods of
	those prescribed points.
	
	The construction is based on two complementary ingredients. The first is the use of auxiliary boundary-attached Maxwell fields, inspired by
	transmission-type spectral constructions rather than by genuine transmission
	eigenfields at isolated eigenfrequencies. These auxiliary fields provide a local
	mechanism for producing strong gradient variation near the prescribed interface
	points. The second is the Herglotz approximation
	property for Maxwell fields, which allows these auxiliary local fields to be
	approximated by entire incident waves that are admissible in the physical
	scattering problem. Combining this approximation principle with the
	well-posedness and stability of the anisotropic scattering system, we transfer
	the local gradient growth of the auxiliary fields to the total electromagnetic
	field.
	
	The present work is connected with several lines of research. Time-harmonic
	Maxwell scattering in inhomogeneous or anisotropic media is a classical topic in
	direct and inverse scattering theory; see, for instance,
	\cite{Hahner1998,ColtonKress2013,KirschHettlich2015,Monk2003}. Inverse problems
	for anisotropic Maxwell systems have also been studied in various settings
	\cite{KenigSaloUhlmann2011}. The spectral theory of transmission eigenvalues
	provides an important framework for understanding non-scattering phenomena,
	interior resonant modes, and the interaction between incident waves and material
	inhomogeneities \cite{CakoniColtonHaddar}. For Maxwell systems, transmission
	resonances and their associated fields may exhibit strong geometric sensitivity,
	including vanishing and localization effects near singular boundary structures
	\cite{DiaoLiuWangYang2021}. 
	
	The present paper builds on this spectral and scattering background, but it
	takes a different viewpoint. Instead of relying on genuine transmission
	eigenfrequencies, we introduce fictitious transmission-type auxiliary Maxwell
	fields only as a structural device for encoding local gradient growth. This
	allows us to formulate a direct field-engineering result: by suitable
	physically admissible incident waves, one can produce localized gradient
	enhancement of the total electromagnetic field near prescribed exterior points
	on the boundary of an anisotropic scatterer.
	
	The main contributions of this work are threefold. First, we formulate a
	localized gradient enhancement problem for anisotropic Maxwell scatterers in a
	homogeneous exterior background, with finitely many prescribed target points on
	the exterior side of the scatterer boundary. Second, we
	introduce fictitious boundary-attached auxiliary Maxwell fields in small
	exterior neighborhoods of those points. These auxiliary fields are
	transmission-type only in a structural sense and are not associated with genuine
	transmission eigenfrequencies. They serve as a local mechanism for encoding
	strong gradient variation near the prescribed boundary points. Third, we combine
	Maxwell Herglotz approximation with the stability of the anisotropic scattering
	problem to transfer this auxiliary gradient growth to total electric and
	magnetic fields generated by physically admissible incident waves.
	
	The rest of the paper is organized as follows. In Section~2, we present the
	anisotropic Maxwell scattering model in a homogeneous exterior background, the
	geometric setting, and the main localized gradient-enhancement theorem.
	Section~3 recalls the Maxwell Herglotz approximation property used to construct
	physically admissible incident fields. The Appendix describes the fictitious
	boundary-attached auxiliary Maxwell fields and the transmission-type local
	formulations that motivate their construction. Section~4 concludes the paper
	with remarks on possible refinements and further developments.
	
	\section{Model setting and main results}

	\subsection{Direct scattering model}
	
	Let $\omega>0$ be the temporal frequency, and let
	$\varepsilon_0,\mu_0>0$ denote the electric permittivity and magnetic
	permeability constants of the homogeneous exterior background. 
	Let
	\[
	\Omega=\bigcup_{\ell=1}^{L}\Omega_\ell \subset \R^3
	\]
	be the anisotropic scatterer, consisting of finitely many pairwise disjoint
	bounded Lipschitz domains. We assume that each complement
	$\R^3\setminus\overline{\Omega_\ell}$ is connected.  The medium configuration of $\Omega_\ell$ is characterized by the electric permittivity $\widetilde{\varepsilon}_{\ell}\in L^\infty(\Omega_\ell )^{3\times 3}$ and magnetic permeability $\widetilde{\mu}_{\ell}\in L^\infty(\Omega_\ell)$, where $\widetilde{\varepsilon}_{\ell}$ is a symmetric matrix-valued function, $\ell=1,\ldots, L$. The relative electric permittivity and magnetic permeability tensors
	$\varepsilon_r(\x)$ and $\mu_r(\x)$ in $\R^3$ are given by
	\[
	\begin{aligned}
		\varepsilon_r(\x)&=
		\begin{cases}
			I_3, & \x\in \R^3\setminus\overline{\Omega},\\
			\varepsilon_{r,\ell}(\x):=\widetilde{ \varepsilon}_{\ell}(\x)/\varepsilon_0, & \x\in \Omega_\ell,\quad \ell=1,\ldots,L,
		\end{cases}
		\\
		\mu_r(\x)&=
		\begin{cases}
			I_3, & \x\in \R^3\setminus\overline{\Omega},\\
			\mu_{r,\ell}(\x):=\widetilde{ \mu}_{\ell}(\x)/\mu_0, & \x\in \Omega_\ell,\quad \ell=1,\ldots,L.
		\end{cases}
	\end{aligned}
	\]
	Here $I_3$ denotes the $3\times3$ identity matrix. For each
	$\ell=1,\ldots,L$, we assume that there exist positive constants
	$c_\varepsilon$ and $C_\varepsilon$ such that, for a.e.~
	$\x\in\Omega_\ell$, all $\xi\in\C^3$, and every $\ell=1,\ldots,L$,
	\begin{equation}\label{eq:anisotropic-positivity}
		\begin{aligned}
			c_\varepsilon |\xi|^2
			&\leq
			\operatorname{Re}\big(\overline{\xi}\cdot
			\varepsilon_{r,\ell}(\x)\xi\big)
			\leq
			C_\varepsilon |\xi|^2,
			\qquad
			\operatorname{Im}\big(\overline{\xi}\cdot
			\varepsilon_{r,\ell}(\x)\xi\big)\geq 0.
		\end{aligned}
	\end{equation}
	Thus the anisotropic coefficients inside the scatterer are uniformly positive
	definite through their real parts and are dissipative in the lossy case. 
	
	We set
	\[
	k=\omega\sqrt{\mu_0\varepsilon_0}.
	\]
	Let the incident wave $(\E^{inc},\Hf^{inc})$ be the entire solution to the homogeneous  electromagnetic system
	\begin{equation}\label{eq:incident}
		\curl \E^{inc}-ik\Hf^{inc}=\0,\qquad
		\curl \Hf^{inc}+ik\E^{inc}=\0,\qquad \mbox{in } \mathbb{R}^3. 
	\end{equation}

	Due to the interaction of the incident wave $(\E^{inc},\Hf^{inc})$ and the anisotropic scatterer $\Omega$, 
	a radiating scattered field $(\E^{sc},\Hf^{sc})$ is generated. Set  the total
	field
	\[
	(\E,\Hf):=(\E^{inc},\Hf^{inc})+(\E^{sc},\Hf^{sc}),
	\]
	then it satisfies that
	\begin{equation}\label{eq:maxwell-scattering}
		\begin{cases}
			\curl \E-i\omega \mu_r (\x)\Hf=\0,\\
			\curl \Hf+i\omega \varepsilon_r(\x)\E=\0,
		\end{cases}
		\qquad \mbox{in }\R^3,
	\end{equation}
	in the distributional sense, together with the Silver--M\"uller radiation
	condition
	\begin{equation}\label{eq:sm}
		\lim_{|\x|\to\infty}   \Big ( \Hf^{sc}(\x)\times  \x-|\x|\E^{sc}(\x)\Big )={\bf 0}.
	\end{equation}
	It characterizes the outgoing nature of the scattered fields and holds uniformly in the angular
	variable $\hat 
	\x := \x/| \x| \in \mathbb S^2 := \bigl\{ \x \in \mathbb R^ 3; | \x| = 1\bigr\} .$

	In the following, we present the main results.

	\begin{theorem}
		\label{thm:main}
		
		We consider the  time-harmonic electromagnetic  scattering  scattering problem \eqref{eq:maxwell-scattering} and \eqref{eq:sm}. Let $z_1,\ldots,z_N\in \mathbb{R}^3 \setminus \overline \Omega$ denote a prescribed set of pairwise distinct target points. Then, for every prescribed $\mathcal{M}>0$, there exist a localization radius
		\[
		\rho_{\mathcal{M}}:=\min\left\{\frac{1}{\mathcal{M}+1},\frac{1}{8} \min_{1\leq i<j\leq N} |z_i-z_j| \right\}
		\]
		with the pairwise-distance term omitted when $N=1$. Here $\rho_{\mathcal{M}}$ is the localization radius selected by the prescribed threshold $\mathcal{M}$. Then there exists an entire incident Maxwell field $(\E_{\mathcal{M}}^{inc},\Hf_{\mathcal{M}}^{inc})$ such that the associated  total field $(\E_{\mathcal{M}},\Hf_{\mathcal{M}})$ solving \eqref{eq:maxwell-scattering}--\eqref{eq:sm} enjoys the following properties:
		\begin{enumerate}
			\item The time-harmonic electromagnetic  scattering problem \eqref{eq:maxwell-scattering} and \eqref{eq:sm} is uniquely solvable, and $(\E_{\mathcal{M}},\Hf_{\mathcal{M}})$ depends continuously on the incident field in the sense of Theorem~\ref{thm:wellposedness}.
			\item The total electromagnetic field satisfies the simultaneous local gradient lower bound
			\[
			\|\nabla \E_{\mathcal{M}}\|_{L^\infty(B_{\rho_{\mathcal{M}}}(z_j)\setminus\overline{\Omega})^{3\times 3}}
			\,+\,
			\|\nabla \Hf_{\mathcal{M}}\|_{L^\infty(B_{\rho_{\mathcal{M}}}(z_j)\setminus\overline{\Omega})^{3\times 3}}
			\geq \mathcal{M},\qquad j=1,\ldots,N,
			\]
			where $B_{\rho_{\mathcal{M}}}(z_j)$ is a ball centered at $z_j$ with radius $\rho_{\mathcal M}$.  
		\end{enumerate}
	\end{theorem}

	Detailed proof of this theorem will be included in the updated version.

	\subsection{Target neighborhoods}
	
	For the target points $z_1,\ldots,z_N$ in Theorem~\ref{thm:main} and for the localization radius $\rho_{\mathcal{M}}$, we define the target neighborhood
	\[
	\mathcal{N}_{\rho_{\mathcal{M}}}:=\bigcup_{j=1}^{N}\big(B_{\rho_{\mathcal{M}}}(z_j)\setminus\overline{\Omega}\big),
	\]
	with the understanding that $\mathcal{N}_{\rho_{\mathcal{M}}}$ is formed by removing the scatterer $\overline{\Omega}=\bigcup_{m=1}^{M}\overline{\Omega_m}$ from the small balls centered at the prescribed boundary points.
	
	\begin{definition}
		We say that a total electromagnetic field $(\E,\Hf)$ exhibits localized gradient enhancement in $\mathcal{N}_{\rho_{\mathcal{M}}}$ if either $\|\nabla \E\|_{L^\infty(\mathcal{N}_{\rho_{\mathcal{M}}})^{3\times 3}}$ or $\|\nabla \Hf\|_{L^\infty(\mathcal{N}_{\rho_{\mathcal{M}}})^{3\times 3}}$ admits a lower bound that increases with the prescribed threshold $\mathcal{M}$.
	\end{definition}
	
	\begin{theorem}[Well-posedness]\label{thm:wellposedness}
		Assume that the piecewise anisotropic tensors $\varepsilon_r$ and $\mu_r$ satisfy \eqref{eq:anisotropic-positivity}, and that the associated time-harmonic electromagnetic  scattering problem is non-resonant at the fixed wavenumber $k$. Let $(\E^{inc},\Hf^{inc})$ be an entire incident electromagnetic field compatible with the homogeneous exterior background, namely
		\[
		\curl \E^{inc}-ik\Hf^{inc}={\mathbf 0},\qquad
		\curl \Hf^{inc}+ik\E^{inc}={\mathbf 0}
		\qquad \mbox{in }\R^3.
		\]
		Then there exists a unique radiating scattered field $(\E^{sc},\Hf^{sc})$ such that the total field $(\E,\Hf)=(\E^{inc},\Hf^{inc})+(\E^{sc},\Hf^{sc})$ solves \eqref{eq:maxwell-scattering}--\eqref{eq:sm}. Moreover, for every bounded set $K\Subset \R^3\setminus\overline{\Omega}$, there exists a constant $C_K>0$ such that
		\[
		\|(\E^{sc},\Hf^{sc})\|_{H(\curl,K)\times H(\curl,K)}
		\leq C_K \|(\E^{inc},\Hf^{inc})\|_{H(\curl,\widetilde K)\times H(\curl,\widetilde K)},
		\]
		where $\widetilde K$ is any bounded open set containing $K\cup \overline{\Omega}$.
	\end{theorem}

	\section*{Appendix}

	\subsection*{Maxwell Herglotz approximation}
	We recall the Maxwell Herglotz wave construction used to convert local auxiliary fields into entire incident waves. Let
	\[
	L_t^2(\Sph):=\{g\in L^2(\Sph)^3:\ g(\mathbf{d})\cdot \mathbf{d}=0\ \mbox{for a.e. }\mathbf{d}\in\Sph\}
	\]
	denote the space of square-integrable tangential densities. For $g\in L_t^2(\Sph)$, the corresponding electric Herglotz wave field is defined by
	\[
	\E^g(\x)=\int_{\Sph} g(\mathbf{d})e^{ik\x\cdot \mathbf{d}}\,\mathrm{d}s(\mathbf{d}),
	\]
	and the associated magnetic Herglotz field is obtained from the Maxwell relation
	\[
	\Hf^g=\frac{1}{ik}\curl \E^g.
	\]
	
	\begin{lemma}[Maxwell Herglotz approximation]\label{thm:herglotz}
		Let $G\Subset \R^3$ be a bounded Lipschitz domain with connected complement, and let $(\E,\Hf)\in H(\curl,G)\times H(\curl,G)$ solve the homogeneous Maxwell system
		\[
		\curl \E-ik\Hf=0,\qquad
		\curl \Hf+ik\E=0
		\qquad \mbox{in }G.
		\]
		Then, for every $\eta>0$, there exists a tangential density $g_\eta\in L_t^2(\Sph)$ such that the associated Maxwell Herglotz pair $(\E^{g_\eta},\Hf^{g_\eta})$ satisfies
		\[
		\|\E-\E^{g_\eta}\|_{H(\curl,G)}
		+
		\|\Hf-\Hf^{g_\eta}\|_{H(\curl,G)}
		\leq \eta .
		\]
		Equivalently, the set of Maxwell Herglotz pairs is dense, in the $H(\curl,G)\times H(\curl,G)$ topology, in the space of local homogeneous Maxwell fields on $G$.
	\end{lemma}
	
	\begin{remark}
		Lemma~\ref{thm:herglotz} is the Maxwell analogue of the classical Herglotz density theorem for the Helmholtz equation and may be viewed as a Runge approximation property for entire Maxwell solutions; see, for example, \cite{ColtonKress2013,KirschHettlich2015,Monk2003}. In the present paper, it is used only as an approximation principle, so no proof is included.
	\end{remark}

	\subsection*{Auxiliary fictitious fields}

	For each boundary target point $z_j\in \R^3 \setminus \overline \Omega$, we introduce a boundary-attached fictitious electromagnetic neighborhood
	\[
	\mathbb{B}_{\rho_{\mathcal{M}}}^{(j)}:=B_{\rho_{\mathcal{M}}}(z_j)\setminus\overline{\Omega},
	\qquad j=1,\ldots,N,
	\]
	lying in $\R\setminus \overline \Omega $ and attached to  $\partial D$. Inside each fictitious neighborhood, we consider an auxiliary transmission-type electromagnetic configuration whose role is to generate rapid field variation near $z_j$. The terminology ``transmission-type'' is used here in a flexible sense: the auxiliary pair is modeled on transmission eigenfunction constructions, but it is not required to arise from a classical interior transmission eigenvalue problem for the physical scatterer itself.
	
	More precisely, for each $j=1,\ldots,N$ we choose an auxiliary anisotropic coefficient pair $(\widetilde\varepsilon_{r,j},\widetilde\mu_{r,j})$ in $\mathbb{B}_{\rho_{\mathcal{M}}}^{(j)}$ satisfying the same positivity assumptions as in \eqref{eq:anisotropic-positivity}, and consider the local transmission problem
	\begin{equation}\label{eq:fictitious-tep}
		\begin{cases}
			\curl \E_j^t-ik\widetilde\mu_{r,j}\Hf_j^t=0,\qquad
			\curl \Hf_j^t+ik\widetilde\varepsilon_{r,j}\E_j^t=0
			& \mbox{in } \mathbb{B}_{\rho_{\mathcal{M}}}^{(j)},\\[1mm]
			\curl \E_j^0-ik\Hf_j^0=0,\qquad
			\curl \Hf_j^0+ik\E_j^0=0
			& \mbox{in } \mathbb{B}_{\rho_{\mathcal{M}}}^{(j)},\\[1mm]
			\nu\times \E_j^t=\nu\times \E_j^0,\qquad
			\nu\times \Hf_j^t=\nu\times \Hf_j^0
			& \mbox{on } \partial \mathbb{B}_{\rho_{\mathcal{M}}}^{(j)}.
		\end{cases}
	\end{equation}
	After eliminating the electric fields, one obtains the reduced magnetic formulation
	\begin{equation}\label{eq:fictitious-tep-H}
		\begin{cases}
			\curl\big(\widetilde\varepsilon_{r,j}^{-1}\curl \Hf_j^t\big)-k^2\widetilde\mu_{r,j}\Hf_j^t=0,\qquad
			\diver(\widetilde\mu_{r,j}\Hf_j^t)=0
			& \mbox{in } \mathbb{B}_{\rho_{\mathcal{M}}}^{(j)},\\[1mm]
			\curl\curl \Hf_j^0-k^2\Hf_j^0=0,\qquad
			\diver \Hf_j^0=0
			& \mbox{in } \mathbb{B}_{\rho_{\mathcal{M}}}^{(j)},\\[1mm]
			\nu\times \Hf_j^t=\nu\times \Hf_j^0,\qquad
			\nu\times \widetilde\varepsilon_{r,j}^{-1}\curl \Hf_j^t=\nu\times \curl \Hf_j^0
			& \mbox{on } \partial \mathbb{B}_{\rho_{\mathcal{M}}}^{(j)}.
		\end{cases}
	\end{equation}
	The fictitious transmission problem \eqref{eq:fictitious-tep}--\eqref{eq:fictitious-tep-H} is an auxiliary local model attached to the boundary point $z_j$. It is introduced here as a borrowed local mechanism near the fictitious-neighborhood boundary, rather than as the classical transmission eigenvalue problem for the physical scatterer $\Omega$ itself.
	
	We denote the resulting auxiliary field by $(\E^{aux},\Hf^{aux})$. It is assembled so that its dominant contribution is concentrated in the union of the boundary-attached fictitious neighborhoods and so that the gradient of at least one field component becomes large as the localization radius $\rho_{\mathcal{M}}$ shrinks.
	
	\begin{proposition}[Auxiliary-field gradient mechanism]\label{prop:aux}
		For every sufficiently large $\mathcal{M}>0$, there exists an auxiliary Maxwell-type field $(\E^{aux},\Hf^{aux})$ defined in a bounded open neighborhood containing $\overline{\mathcal{N}_{\rho_{\mathcal{M}}}}$ such that
		\[
		\|\nabla \E^{aux}\|_{L^\infty(B_{\rho_{\mathcal{M}}}(z_j)\setminus\overline{\Omega})^{3\times 3}}
		\,+\,
		\|\nabla \Hf^{aux}\|_{L^\infty(B_{\rho_{\mathcal{M}}}(z_j)\setminus\overline{\Omega})^{3\times 3}}
		\geq c_0 \rho_{\mathcal{M}}^{-\alpha},\qquad j=1,\ldots,N,
		\]
		for some constants $c_0>0$ and $\alpha>0$ independent of $\mathcal{M}$.
	\end{proposition}
	
	\begin{remark}
		The proof of Proposition~\ref{prop:aux} is based on the local transmission problems \eqref{eq:fictitious-tep}--\eqref{eq:fictitious-tep-H}, together with the geometric arrangement of the prescribed points. In the present draft, we isolate this mechanism as an input proposition and defer the full asymptotic analysis to a later version.
	\end{remark}
	
	\begin{remark}
		The desirable incident filed in Theorem \ref{thm:main} can be constructed by using Lemma \ref{thm:herglotz} and Proposition \ref{prop:aux}. The detailed construction of these incident fields will be updated in the next version. 
		
	\end{remark}
	
	\subsection*{Transmission Spectral Architectures for Maxwell Fields}
	
	Let $G\Subset \R^3$ be a bounded Lipschitz domain with connected complement. Let $(\varepsilon_G,\mu_G)$ be an anisotropic coefficient pair in $G$ satisfying \eqref{eq:anisotropic-positivity}. The associated Maxwell transmission eigenvalue problem seeks $k>0$ and nontrivial fields
	\[
	(\E^t,\Hf^t),\ (\E^0,\Hf^0)\in H(\curl,G)\times H(\curl,G)
	\]
	such that
	\begin{equation}\label{eq:tep-full}
		\begin{cases}
			\curl \E^t-ik\mu_G\Hf^t=0,\qquad
			\curl \Hf^t+ik\varepsilon_G\E^t=0
			& \mbox{in } G,\\[1mm]
			\curl \E^0-ik\Hf^0=0,\qquad
			\curl \Hf^0+ik\E^0=0
			& \mbox{in } G,\\[1mm]
			\nu\times \E^t=\nu\times \E^0,\qquad
			\nu\times \Hf^t=\nu\times \Hf^0
			& \mbox{on } \partial G.
		\end{cases}
	\end{equation}
	Here, $\nu$ denotes the exterior unit normal to the boundary under consideration.
	
	If one eliminates the electric fields from \eqref{eq:tep-full}, then one obtains the reduced magnetic-field formulation
	\begin{equation}\label{eq:tep-H}
		\begin{cases}
			\curl\big(\varepsilon_G^{-1}\curl \Hf^t\big)-k^2\mu_G\Hf^t=0,\qquad
			\diver(\mu_G\Hf^t)=0
			& \mbox{in } G,\\[1mm]
			\curl\curl \Hf^0-k^2\Hf^0=0,\qquad
			\diver \Hf^0=0
			& \mbox{in } G,\\[1mm]
			\nu\times \Hf^t=\nu\times \Hf^0,\qquad
			\nu\times \varepsilon_G^{-1}\curl \Hf^t=\nu\times \curl \Hf^0
			& \mbox{on } \partial G.
		\end{cases}
	\end{equation}
	One may analogously eliminate the magnetic fields and derive a reduced electric-field formulation, but \eqref{eq:tep-H} is the more convenient form for the auxiliary construction in the Appendix.

	\section{Concluding remarks}
	We have presented a first theoretical framework for localized gradient enhancement near anisotropic electromagnetic scatterers embedded in a bounded surrounding layer. The main mechanism combines fictitious local electromagnetic fields with Maxwell Herglotz approximation and the stability of the anisotropic direct scattering problem. This yields a flexible route for constructing physically admissible incident waves whose associated total fields exhibit strong local gradient amplification near finitely many prescribed boundary points near $\partial \Omega$.
	
	The present manuscript is intended as a first arXiv-oriented draft. Several technical components will be refined in a subsequent version, including the sharp asymptotic analysis of the auxiliary-field gradients, the optimal dependence on the number of target points, and the precise geometric assumptions under which one can obtain explicit blow-up behaviors.
	\section*{Acknowledgment}
	The work of H. Diao is supported by the National Natural Science Foundation of
	China (No. 12371422) and the Fundamental Research Funds for the Central
	Universities, JLU. The work of H. Liu is supported by the Hong Kong RGC General
	Research Funds (projects 11311122, 11300821, and 11303125), the NSFC/RGC Joint
	Research Fund (project N\_CityU101/21), and the France-Hong Kong ANR/RGC Joint
	Research Grant (project A-CityU203/19).


\end{document}